\documentclass[a4paper,11pt,reqno]{amsart} 
\usepackage{amsmath,amssymb,amsthm}
\usepackage{graphpap,latexsym,epsf}
\usepackage{epsf}
\usepackage{color,graphics}
\usepackage[numbers,sort&compress,longnamesfirst]{natbib}

\theoremstyle{plain}
\newtheorem{theorem}{Theorem}[section]  
\newtheorem{corollary}[theorem]{Corollary}  
\newtheorem{lemma}[theorem]{Lemma}  

\newtheorem*{theorem*}{Theorem}

\theoremstyle{plain}  
\newtheorem{remark}[theorem]{Remark}

\newtheoremstyle{citing}
  {3pt}
  {3pt}
  {\itshape}
  {}
  {\bfseries}
  {.}
  {.5em}
  {\thmnote{#3}}

\theoremstyle{citing}

\numberwithin{equation}{section}

\newcommand{\dist}{\text{dist}}

\newcommand{\rand}{\partial}

\newcommand{\intg}{\mathop{\int}\limits}  
  
\newcommand{\supp}{\mbox{supp}}

\newcommand{\di}{\;d}  
  
\newcommand{\nz}{{\mathbb N}}

\newcommand{\rz}{{\mathbb R}}

\newcommand{\eps}{\varepsilon}  
\renewcommand{\phi}{\varphi}

\renewcommand{\div }{{\rm div}\,}
\renewcommand{\a }{\alpha }

\newcommand{\media}{\mkern12mu\hbox{\vrule height4pt depth-3.2pt width5pt} \mkern-16mu\int}

\newcommand{\Di }{\mathcal{D}^{1,2}_{a}(\Omega) }

\newcommand{\n}{\nabla}

\newcommand{\R}{\mathbb{R}} 
 
\newcommand{\N}{\mathbb{N}} 
\newcommand{\caratt}{{\setbox0 =\hbox{$\mathsurround=0pt\chi$}\hbox
  {\raise\dp0 \copy0 }}}

\begin{document}
 
\title[Regularity for degenerate elliptic equations]
{A note on regularity of solutions to degenerate elliptic equations of Caffarelli-Kohn-Nirenberg type}

\author{Veronica Felli}
\author{Matthias Schneider}
\thanks{ V. F. is supported by M.U.R.S.T. under the national project
  ``Variational Methods and Nonlinear Differential Equations''. During
  the preparation of the paper M. S. research was supported by a
  S.I.S.S.A. postdoctoral fellowship.} 
\address{Scuola Internazionale Superiore di Studi Avanzati\\
S.I.S.S.A.\\
Via Beirut 2-4\\
34014 Trieste, Italy}
\email{felli@sissa.it} 
\address{Mathematisches Institut\\
    Im Neuenheimer Feld 288\\
   D-69120 Heidelberg, Germany}
\email{mschneid@mathi.uni-heidelberg.de}

\date{June 25, 2003}  

\begin{abstract}
We establish H\"older continuity of weak solutions to degenerate critical
elliptic equations of Caffarelli-Kohn-Nirenberg type. 
\end{abstract}

\maketitle


\section{Introduction}
\label{s:intro}

Our purpose is to establish H\"older continuity of weak solutions to
\begin{align}
\label{eq:23}
-\div(|x|^{-2a} \nabla u) = \frac{f}{|x|^{bp}}, \text{ in } \Omega \subset \rz^N,   
\end{align}
where $\Omega$ is an open, $N \ge 3$ and $a$, $b$, and $p$ satisfy
\begin{align}
\label{eq:24}
\begin{split}
&-\infty<a<\frac{N-2}{2},\quad a\le b \le a+1\\ 
&p=p(a,b)= \frac{2N}{N-2(1+a-b)}.  
\end{split}
\end{align}
We denote by $\mathcal{D}^{1,2}_{a}(\Omega)$ the closure of $C_c^{\infty}(\Omega)$ with
respect to the norm 
\begin{align*}
\|u\|_{\mathcal{D}^{1,2}_{a}(\Omega)} := \left(\int_\Omega |\nabla u|^2 |x|^{-2a}\right)^{1/2}.   
\end{align*}
For a given weight $\omega$ we denote by $L^p(\Omega,\omega)$ the space of functions $u$ such that
\begin{align*}
\|u\|^p_{L^p(\Omega,\omega)} := \int_\Omega |u|^p \omega(x) < \infty.   
\end{align*}
The space $H^1_a(\Omega)$ is defined to be the closure of $C^\infty(\bar \Omega)$ with respect to 
\begin{align*}
\|u\|^2_{H^1_a(\Omega)} := \int_\Omega |x|^{-2a} \left(|\nabla u|^2+|u|^2\right).    
\end{align*}
Our interest in these problems arose because of their relation to nonlinear, degenerate
elliptic equations stemming from the family of Caffarelli-Kohn-Nirenberg inequalities
\cite{CKN}: if $a$, $b$, and $p$ satisfy (\ref{eq:24}) then we have for all $u \in
\mathcal{D}^{1,2}_{a}(\rz^N)$   
\begin{align}
\label{eq:25}
\left(\int |u|^p |x|^{-bp}\right)^{1/p} \le {\mathcal C}_{a,b,N} \left(\int |\nabla u|^2 |x|^{-2a}\right)^{1/2}.  
\end{align}
For best constants and existence of minimizers in (\ref{eq:25}) we refer to
\cite{CatrinaWang}. Due to its characterization any minimizer $u \in
\mathcal{D}^{1,2}_{a}(\rz^N)$, if it exists, is a weak
solution to 
\begin{align}
\label{eq:26}
-\div(|x|^{-2a} \nabla u) = \frac{K(x) |u|^{p-2}u}{|x|^{bp}} \text{ in } \rz^N,  
\end{align}
where $K(x) \equiv {\rm const}$ for some appropriate constant. 
The exponent $p=p(a,b,N)$ is the critical exponent in (\ref{eq:25}) and shares many features with the critical Sobolev exponent, e.g. (\ref{eq:26}) possesses for $K(x) \equiv K(0)$ a dilation symmetry,
which gives rise to a noncompact manifold of weak radial solutions for $K(0)>0$. 
In order to study problem (\ref{eq:26}) for non-constant functions $K$ using for instance a degree argument, H\"older estimates
for weak solutions of (\ref{eq:23}) are an important tool (see \cite{SF_bup}).\\
Regularity properties of weak solution to degenerate elliptic problems 
with more general weighted operators of the form $\div(\omega(x)\nabla(\cdot))$
are studied in \cite{CiccoVivaldi96,FabesKenigSera82,potential} (see also the references mentioned there). The classes of weights $\omega$ treated there include the class $(QC)$ of weights 
\begin{align*}
\omega(x)= |{\det} T'|^{1-2/N},
\end{align*} 
where $T:\rz^N \to \rz^N$ is quasi-conformal (see \cite{FabesKenigSera82,potential} for a definition). In fact our weights $|\cdot|^{-2a}$ are associated with quasi-conformal transformations $T_a(x):= x |x|^{-2a/(N-2)}$. The right-hand sides studied in \cite{CiccoVivaldi96,FabesKenigSera82,potential} are either zero or in divergence form, e.g. H\"older continuity of weak solutions to
\begin{align*}
-\div(|x|^{-2a} \nabla u) = \div(F) \text{ in } \Omega
\end{align*}  
is established in \cite{CiccoVivaldi96} assuming $|F| |x|^{2a}\in  L^p(\Omega,|x|^{-2a})$ for some $p>\max(N-2a,N,2)$. We derive H\"older estimates for weak solutions to (\ref{eq:23}) in terms of $f$, because a sharp relation of the integrability of $f$ and its representation in divergence form $F$ in the various weighted spaces is not obvious. We compare weak solutions of (\ref{eq:23}) with $\mu_a$-harmonic functions, which are by definition weak solutions of
\begin{align*}
-\div(|x|^{-2a} \nabla u) = 0 \text{ in } \Omega,
\end{align*}  
where H\"older regularity is known (see for instance \cite{potential}) and prove 
\begin{theorem}\label{t:holder_estimate}
Suppose $\Omega \subset \rz^N$ is a bounded domain and $u \in H^1_{a}(\Omega)$ weakly solves
(\ref{eq:23}), that is
\[
\int_{\Omega}|x|^{-2a}\n u\cdot\n \varphi\,dx=
\int_{\Omega}|x|^{-bp}f\varphi\,dx\quad\forall\,\varphi\in
H^1_{0,a}(\Omega).
\]
Assume $a$, $b$ and $p$ satisfy (\ref{eq:24}), $b<a+1$, and 
$f\in L^{s}(\Omega,|x|^{-bp})$ for some $s>p/(p-2)$.
Then $u\in C^{0,\a}$ for any $\alpha \in (0,1)$ satisfying
\[\a<\min(\a_h,1) \text{ and }
\a<\begin{cases}
\left(\frac{N-2}2-a\right)\left(p-2-\frac{p}{s}\right) &\text{if }b\ge 0\\
\frac{N}{p}\left(p-2-\frac{p}{s}\right) &\text{if }b< 0
\end{cases}
,\] 
where $\a_h$ is the regularity of
$\mu_a$-harmonic functions given in Theorem \ref{t:mu-harm-funct} below. 
Moreover, for any $\Omega'\Subset \Omega$ there is a constant $C=C(N,a,\a,\Omega,\dist(\Omega',\Omega),s)$ 
such that 
\[
\sup_{\Omega'}|u|+\sup_{\substack{x,y\in\Omega'\\x\not=y}}\frac{|u(x)-u(y)|}{|x-y|^{\a}}\leq
C\,\Big\{\|u\|_{L^2(\Omega,\di\mu_a)}+\|f\|_{L^{s}(\Omega,|x|^{-bp})}\Big\}
\]
\end{theorem}
For the nonlinear problem (\ref{eq:26}) we use a De Giorgi-Moser type iteration procedure
as in \cite{BrezisKato} and obtain
\begin{theorem}\label{t:reg}
Let $a$, $b$ and $p$ satisfy (\ref{eq:24}) and $u\in D_a^{1,2}(\R^N)$ be a weak solution to
\begin{equation}
\label{eq:6bk}
-\div(|x|^{-2a}\n u)=K(x)\frac{|u|^{p-2}u}{|x|^{bp}}, \quad x\in\Omega
\end{equation}
where $K\in L^{\infty}(\Omega)$. Then $u\in L^s_{\rm loc}(\Omega,|x|^{-bp})$ for
any $s\in[p,+\infty[$. Moreover, $u$ is H\"older continuous with H\"older exponent given in
Theorem \ref{t:holder_estimate}.
\end{theorem}
 
\begin{remark}
While completing this note we learned that in \cite{ColoPeral} 
weighted $q$-Laplacian equations of the form
$$
-\div(|x|^{-qa} |\nabla u|^{q-2}\nabla u)=g
$$
are studied. Under assumption (\ref{eq:24}) H\"older regularity of 
weak solutions to equation
\begin{align*}
-\div(|x|^{-qa} \nabla u) = \frac{f}{|x|^{bp}}, \text{ in } \Omega \subset \rz^N,   
\end{align*}
is shown if $a=b$, $a>-1$, and 
$f\in L^{s}(\Omega,|x|^{-bp})$ for some $s>p/(p-2)$. Theorem
\ref{t:holder_estimate} extends  
this result to the full range for $a$ and $b$ in the case $q=2$.    
\end{remark}

\section{Preliminaries}
We collect some properties of the weighted measure $\mu_a:=|x|^{-2a}\, dx$ and $\mu_a$-harmonic functions. 
We refer to \cite{FabesKenigSera82,potential} for the proofs.
\begin{itemize}
\item The measure $\mu_a$ satisfies the doubling property, i.e. for every $\tau \in (0,1)$              
there exists a constant $C_{(\ref{eq:3})}(\tau)$ such that
        \begin{align}
        \label{eq:3}
        \mu_a(B(x,r)) \le C_{(\ref{eq:3})}(\tau) \mu_a(B(x,\tau r))  
        \end{align}
\item A Poincar\'e-type inequality holds, i.e. 
there is a positive constant $C_{(\ref{eq:8})}$ such that any $u
\in D_a^{1,2}(\rz^N)$ satisfies  
\begin{align}
\label{eq:8}
\media_{B_r(x)}|u-u_{x,r}|^2 \di\mu_a \leq
C_{(\ref{eq:8})} \, r^2 \media_{B_r(x)} |\n u|^2 \di\mu_a,
\end{align}
where $u_{x,r}$ denotes the weighted mean-value 
\[
u_{x,r}:=\media_{B_r(x)}u \di \mu_a = \frac{1}{\mu_a(B_r(x))}\int_{B_r(x)}u(x)\di\mu_a.
\]
\end{itemize}
Concerning $\mu_a$- harmonic functions we have the following results. 
\begin{theorem}[Thm. 3.34 in \cite{potential}(p. 65), Thm. 6.6 in \cite{potential}(p. 111)]
\label{t:mu-harm-funct} 
$ $\\
There are constants $C_{(\ref{eq:6})}(N,a)$ and $\alpha_h=\alpha_h(N,a) \in (0,1)$ such that 
if $u$ is $\mu_a$-harmonic in $B_r(x_0) \subset \rz^N$ and $0<\rho<r$ then
\begin{align}
\label{eq:6}
\text{ess-sup}_{B(x_0,\frac{r}{2})} |u| \le C_{(\ref{eq:6})} \media_{B(x_0,r)} |u|^2 \di\mu_a,     
\end{align}
\begin{align}
\label{eq:7}
osc(u,B_\rho(x_0)) &\le 2^{\alpha_h} \big(\frac{\rho}{r}\big)^\alpha osc(u,B_r(x_0)).  
\end{align}
Consequently, $\mu_a$-harmonic functions are H\"older continuous.
\end{theorem}
We will call a function $u \in D_{a,loc}^{1,2}(\rz^N)$ weakly super $\mu_a$-harmonic in $\Omega$, if
 for all nonnegative $\varphi \in C^\infty_c(\Omega)$ we have
\begin{align}
\label{eq:11}
\int_\Omega |x|^{-2a} \nabla u \nabla \varphi \ge 0.  
\end{align}
 
\begin{theorem}[Thm 3.51 in \cite{potential}(p. 70)]
\label{thm-weak-harnack}
There exist positive constants $s=s(N,a)$ and 
$C_{(\ref{eq:12})}=C_{(\ref{eq:12})}(N,a)$ such that
if $u$ is nonnegative and weakly super $\mu_a$-harmonic in $\Omega$ and $B_{2r}(x_0) \subset \Omega$ we have
\begin{align}
\label{eq:12}
\text{ess inf}_{B_{\frac{r}{2}}(x_0)} u\ge C_{(\ref{eq:12})} \big(\media_{B_{r}(x_0)} u^s \di
\mu_a\big)^{\frac{1}{s}}.    
\end{align}
\end{theorem}
We use the two theorems above to derive 
\begin{lemma}
\label{reg:lem:2}
For any ball $B_r(x_0)$ there is a constant $C_{(\ref{eq:4})}(B_r(x_0))$ such that any $\mu_a$-harmonic
function $u$ in $B_r(x_0)$ satisfies
\begin{align}
\label{eq:4}
\media_{B_\rho(x_0)} |\nabla u|^2\di\mu_a  \le C_{(\ref{eq:4})}
\big(\frac{\rho}{r}\big)^{2\alpha_h -2} \media_{B_r(x_0)} |\nabla u|^2 \di\mu_a,  
\end{align}
where $\a_h \in (0,1)$ is given in Theorem \ref{t:mu-harm-funct}.
\end{lemma}
\begin{proof}
To prove the claim we may assume that $0<\rho<(1/4)r$ and that $u$ has mean-value zero in $B_r(x_0)$.
We take a cut-off function $\xi \in
C_c^\infty(B_{2\rho}(x_0))$ such that $\xi \equiv 1$ in $B_\rho(x_0)$, $0\le \xi \le 1$,
$\|\nabla \xi\|_\infty \le 2 \rho^{-1}$ and
define $\phi:= \xi^2 (u-u(x_0))$. Testing with $\phi$ and using H\"older's inequality we get 
\begin{align}
\intg_{B_r(x_0)} |\nabla u|^2 \xi^2 \di \mu_a \le 
\intg_{B_r(x_0)} |\nabla \xi|^2 (u-u(x_0))^2 \di \mu_a 
\le \|u-u(x_0)\|_{\infty,B_{2\rho}(x_0)}^2 \mu_a\big(B_{2\rho}(x_0)\big) \rho^{-2}.
\label{eq:13}
\end{align}
From (\ref{eq:13}) and (\ref{eq:7}) we infer
\begin{align*}
\media_{B_\rho(x_0)} |\nabla u|^2 \di\mu_a &\le C \big(\frac{\rho}{r}\big)^{2\alpha} 
osc(u, B_{\frac{r}{2}}(x_0))^2  \rho^{-2}\le C \big(\frac{\rho}{r}\big)^{2\alpha}  \rho^{-2} 
\media_{B_r(x_0)} |u|^2 \di\mu_a  
\end{align*}
Finally, since $u$ has mean-value zero in $B_r(x_0)$ the Poincar\'e inequality (\ref{eq:8})
yields the claim.
\end{proof}

\section{Growth of local integrals}
We give a weighted version of the Campanato-Morrey characterization of H\"older continuous functions.
\begin{theorem}\label{t:avera}
Suppose $\Omega \subset \rz^N$ is a bounded domain and $u\in L^2(\Omega,d\mu_a)$ satisfies
\begin{equation}
\label{eq:ass1}
\media_{B_r(x)}|u(y)-u_{x,r}|^2\di\mu_a\leq
M^2 r^{2\a}\quad\text{for any}\ B_r(x)\subset\Omega
\end{equation}
and some $\a\in(0,1)$. Then $u\in C^{0,\a}(\Omega)$ and for any
$\Omega'\Subset\Omega$ there holds
\[
\sup_{\Omega'}|u|+\sup_{\substack{x,y\in\Omega'\\x\not=y}}\frac{|u(x)-u(y)|}{|x-y|^{\a}}\leq
C\Big\{M+\|u\|_{L^2(\Omega,|x|^{-2a})}\Big\}
\]
where $C=C(N,a,\a,\Omega,\dist(\Omega',\Omega))$ is a positive constant independent of $u$.
\end{theorem}
\begin{proof}
Denote $R_0=\dist(\Omega',\partial\Omega)$. Using the triangle inequality and integrating in $B_{r_1}(x_0)$ we have for any $x_0\in\Omega'$ and $0<r_1<r_2\leq R_0$ 
\begin{align*}
&\big|u_{x_0,r_1}-u_{x_0,r_2}\big|^2\\
&\qquad \leq
\frac{2}{\mu_a(B_{r_1}(x_0))}\bigg\{\int_{B_{r_1}(x_0)}\big|u(x)-u_{x_0,r_1}\big|^2\di\mu_a
+\int_{B_{r_2}(x_0)}\big|u(x)-u_{x_0,r_2}\big|^2 \di\mu_a\bigg\}.
\end{align*}
Using assumption (\ref{eq:ass1}) we obtain
\begin{equation}\label{eq:diff1}
\big|u_{x_0,r_1}-u_{x_0,r_2}\big|^2\leq
\frac{2M^2}{\mu_a(B_{r_1}(x_0))}\big\{\mu_a(B_{r_1}(x_0))r_1^{2\a}+\mu_a(B_{r_2}(x_0))r_2^{2\a}\big\}.
\end{equation}
For any $R\leq R_0$ we take $r_1=2^{-(i+1)}R$ and $r_2=2^{-i}R$ in
(\ref{eq:diff1}). The doubling property (\ref{eq:3}) then gives
\begin{align*}
&\big|u_{x_0,2^{-(i+1)}R}-u_{x_0,2^{-i}R}\big|
\leq2M^2\Big(1+C_{(\ref{eq:3})}(N,a)2^{2\a}\Big)2^{-2(i+1)\a}R^{2\a}.
\end{align*}
We sum up and get for $h<k$
\begin{align}\label{eq:diff2}
\big|u_{x_0,2^{-h}R}-u_{x_0,2^{-k}R}\big| &\leq \frac{C(N,a,\a)M}{2^{h\a}}R^{\a}.
\end{align}
The above estimates prove that $\{u_{x_0,2^{-i}R}\}_{i\in\N}\subset\R$
is a Cauchy sequence in $\R$, hence it converges to some limit,
denoted as $\hat u(x_0)$. The value of $\hat u(x_0)$ is independent
of $R$, which may be seen by analogous estimates.
Consequently, from (\ref{eq:diff2}) we have that
\begin{equation}\label{eq:diff3}
\big|u_{x_0,r}-\hat u(x_0)\big|\leq C(N,a,\a)Mr^{\a}\quad\forall\,
x_0\in\Omega'.
\end{equation}
By the Lebesgue theorem we infer  
\[
u_{x,r}=\frac{|B_r(x)|}{\int_{B_r(x)}|y|^{-2a}\,dy}\cdot
\frac{\int_{B_r(x)}|y|^{-2a}u(y)\,dy}{\frac1{|B_r(x)|}}\
\mathop{\longrightarrow}\limits_{r\to 0^+}\
u(x),\quad\text{a. e. in}\ \Omega'.
\]
Hence $\hat u=u$ a. e. in $\Omega'$ and (\ref{eq:diff3}) gives
\begin{equation}\label{eq:diff4}
\big|u_{x_0,r}- u(x_0)\big|\leq C(N,a,\a)Mr^{\a}\quad\forall\,
x_0\in\Omega',
\end{equation}
which implies that $u_{x,r}$ converges to $u$ uniformly in
$\Omega'$. Since $x\mapsto u_{x,r}$ is a continuous function, we conclude that $u$ is continuous in $\Omega'$. From (\ref{eq:diff4}) we have 
\[
|u(x)|\leq C(N,a,\a) MR^{\a}+|u_{x,R}|\quad\forall\, x\in\Omega',\
 \forall \,R\leq R_0.
\]
Thus $u$ is bounded in $\Omega'$ with the estimate
\begin{equation}\label{eq:diff5}
\sup_{\Omega'}|u|\leq c(N,a,\a,\Omega,\dist(\Omega',\Omega))\Big\{M+\|u\|_{L^2(\Omega,|x|^{-2a})}\Big\}.
\end{equation}
Let us now prove that $u$ is H\"older continuous. Let $x,y\in\Omega'$
with $|x-y|=R<\frac{R_0}2$. Assume that $|x|<|y|$. Then we have 
\[
|u(x)-u(y)|\leq|u(x)-u_{x,2R}|+|u(y)-u_{y,2R}|+|u_{x,2R}-u_{y,2R}|.
\]
The first two terms are estimated by (\ref{eq:diff4}), whereas for the
last term we have 
\[
|u_{x,2R}-u_{y,2R}|^2\leq
 2\big\{|u_{x,2R}-u(\xi)|^2+|u(\xi)-u_{y,2R}|^2\big\}
\]
and integrating with respect to $\xi$ over $B_{2R}(x)\cap
  B_{2R}(y)\supseteq B_R(x)$ we obtain
\begin{align*}
&|u_{x,2R}-u_{y,2R}|^2\leq\frac{2}{\mu_a(B_{R}(x))}
\Big(M^2\mu_a(B_{2R}(x))2^{2\a}R^{2\a}+M^2\mu_a(B_{2R}(y))2^{2\a}R^{2\a}\Big).
\end{align*}
Since $x$ is closer to $0$ than $y$, we have that $\mu_a(B_{2R}(y))\leq \mu_a(B_{2R}(x))$ and hence 
\[
|u(x)-u(y)|\leq C(N,a,\a)M|x-y|^{\a}.
\]
If $|x-y|>\frac{R_0}2$ we can use estimate (\ref{eq:diff5}) thus
finding
\begin{align*}
|u(x)-u(y)|&\leq 2\sup_{\Omega'}|u|\leq c2^{\a}
\bigg[M+\frac1{R_0^{\a}}\|u\|_{L^2(\Omega,|x|^{-2a})}\bigg]|x-y|^{\a}.
\end{align*}
The proof is thereby complete.
\end{proof}
\begin{corollary}\label{c:avera}
Suppose $\Omega \subset \rz^N$ is a bounded domain and $u\in H^1_{a}(\Omega)$ satisfies
\[
\media_{B_r(x)}|\n u|^2 \di\mu_a \leq
M^2 r^{2\a-2}\quad\text{for any}\ B_r(x)\subset\Omega
\]
and some $\a\in(0,1)$. Then $u\in C^{0,\a}(\Omega)$ and for any
$\Omega'\Subset \Omega$ there holds
\[
\sup_{\Omega'}|u|+\sup_{\substack{x,y\in\Omega'\\x\not=y}}\frac{|u(x)-u(y)|}{|x-y|^{\a}}\leq
c\Big\{M+\|u\|_{L^2(\Omega,|x|^{-2a})}\Big\}
\]
where $c=c(N,a,\a,\Omega,\dist(\Omega',\Omega))>0$.
\end{corollary}
\begin{proof}
The proof follows from Theorem \ref{t:avera} and the Poincar\'e type inequality in (\ref{eq:8}).
\end{proof}

\begin{proof}[{\bf Proof of Theorem \ref{t:holder_estimate}}]
Let $w\in u+H^1_{0,a}(B_r(x_0))$ be the unique
solution to the Dirichlet problem
\begin{align}
\label{eq:14}
\begin{cases}
-\div(|x|^{-2a}\n w)=0&\text{in}\ B_r(x_0)\\
w\big|_{\partial B_r(x_0)}=u.
\end{cases}  
\end{align}
Clearly the function $v=u-w\in  H^1_{0,a}(B_r(x_0))$ weakly
solves
\[
-\div(|x|^{-2a}\n v)=\frac{f}{|x|^{bp}}\quad\text{in } B_r(x_0).
\]
Testing the above equation with $v$ and using H\"older's inequality and (\ref{eq:25}), we
get
\begin{align*}
\int_{B_r(x_0)}|\n v|^2\di\mu_a x
&\leq {\mathcal C}_{a,b,N} 
\bigg(\int_{ B_r(x_0)}|x|^{-bp}|f|^{\frac{p}{p-1}}\bigg)^{\frac{p-1}p}
\bigg(\int_{ B_r(x_0)}|\n v|^2\di \mu_a \bigg)^{\frac{1}{2}}.
\end{align*}
Since $f\in L^{s}(\Omega,|x|^{-bp})$ for some $s>p/(p-2)$
we may use H\"older's inequality with conjugate exponents $s(p-1)/p$ and
\[
\frac{s(p-1)}{s(p-1)-p} = \frac{p-1}{1+(p-2-\frac{p}{s})}
\]
and Lemma \ref{l:veronica} with $\eps=2(p-2-p/s)/p $ to obtain
\begin{align}
\notag
\int_{ B_r(x_0)}|\n v|^2\di\mu_a &\leq 
{\mathcal C}^2_{a,b,N} \bigg(\int_{B_r(x_0)}|x|^{-bp} |f|^{s}\bigg)^{\frac{2}s}
\bigg(\int_{B_r(x_0)} |x|^{-bp} \bigg)^{\frac{2}{p}+\eps}\\ 
\label{eq:16}
&\le C r^{-2+N\eps}\max(r,|x_0|)^{-bp\eps}\mu_a(B_r(x_0))
\left(\int_{B_r(x_0)}{|x|^{-bp}} |f|^{s} \right)^{\frac{2}s}.
\end{align}
From (\ref{eq:4}) we deduce for any $0<\rho\leq r$
\begin{align}
\notag
\media_{B_{\rho}(x_0)}|\n u|^2\di\mu_a &\leq 4\media_{B_{\rho}(x_0)}|\n w|^2\di\mu_a+
4\media_{B_{\rho}(x_0)}|\n v|^2\di\mu_a\\
&\le 4 C_{(\ref{eq:4})} \big(\frac{\rho}{r}\big)^{2\alpha_h -2} \media_{B_r(x_0)} |\nabla w|^2
\di\mu_a + 4 \mu_a(B_\rho(x_0))^{-1} \int_{B_{r}(x_0)}|\n v|^2\di\mu_a
\label{eq:15}
\end{align}
Since $w$ minimizes the Dirichlet integral we may replace $w$ in (\ref{eq:15}) by $u$.
If we further estimate the integral containing $v$ in (\ref{eq:15}) using (\ref{eq:16}) we get    
\begin{align*}
\int_{B_\rho(x_0)} |\n u|^2\di \mu_a \leq
C\bigg(&\mu_a(B_{\rho}(x_0))\mu_a(B_{r}(x_0))^{-1}\bigg(\frac{\rho}r\bigg)^{-2+2\a_h}
\int_{B_r(x_0)}|\n u|^2\di\mu_a \\
&+r^{-2+N\eps} \max(r,|x_0|)^{-bp\eps} \mu_a(B_{r}(x_0) \|f\|^2_{L^{s}(B_r(x_0),|x|^{-bp})}\bigg).
\end{align*}
We estimate the term $\max(r,|x_0|)^{-bp\eps}$ by $r^{-bp\eps}$ if $b\ge 0$ and in the case $b<0$ by a constant $C(\Omega)$. For the rest of the proof we will consider the more interesting situation $b\ge 0$. The case $b<0$ may be treated analogously.\\  
Lemma \ref{reg:lem:1} with $\Phi(\rho):= \int_{B_\rho(x_0)} |\n u|^2\di \mu_a$ 
gives for $0<\rho<r\le r_0:=\dist(x_0, \rand \Omega)$
\begin{align*}
\int_{B_{\rho}(x_0)} |\n u|^2\di\mu_a\leq
C(\alpha)
\bigg(&\frac{\mu_a(B_{\rho}(x_0))}{\mu_a(B_{r}(x_0))}\bigg(\frac{\rho}r\bigg)^{-2+2\a}
\int_{B_r(x_0)}|\n u|^2\di\mu_a\\
&+{\rho}^{-2+(N-bp)\eps}\mu_a(B_{{\rho}}(x_0)\|f\|^2_{L^{s}(\Omega,|x|^{-bp})}\bigg).
\end{align*}
We take a cut-off function $\xi \in
C_c^\infty(B_{r}(x_0))$ such that $\xi \equiv 1$ in $B_{r/2}(x_0)$, $0\le \xi \le 1$,
$\|\nabla \xi\|_\infty \le 2 r^{-1}$ and
define $\phi:= \xi^2 u$. Testing with $\phi$ and using (\ref{eq:25}) and H\"older's inequality we get 
\begin{align*}
\int_{B_r(x_0)} |\nabla u|^2 \xi^2 \di \mu_a 
&\le {\mathcal C}_{a,b,N} \|f\|_{L^{\frac{p}{p-1}}(\Omega,|x|^{-bp})}  
\|\xi^2 u\|_{\mathcal{D}^{1,2}_{a}(\Omega)} \\
&\quad+ \|u \nabla\xi\|_{L^{2}(B_r(x_0),\di\mu_a)}
\|\nabla u \xi\|_{L^{2}(B_r(x_0),\di\mu_a)}.
\end{align*}
We divide by $\|\nabla u \xi\|_{L^{2}(B_r(x_0),\di\mu_a)}$ and obtain
\begin{align*}
\int_{B_r(x_0)} |\nabla u|^2 \xi^2 \di \mu_a 
&\le 
\frac{2{\mathcal C}^2_{a,b,N}\|f\|^2_{L^{\frac{p}{p-1}}(\Omega,|x|^{-bp})} 
\|\xi^2 u\|^2_{\mathcal{D}^{1,2}_{a}(\Omega)}}
{\|\nabla u \xi\|^2_{L^{2}(B_r(x_0),\di\mu_a)}}
+\|u \nabla\xi\|^2_{L^{2}(B_r(x_0),\di\mu_a)}\\
&\le \frac{2{\mathcal C}^2_{a,b,N}\|f\|^2_{L^{\frac{p}{p-1}}(\Omega,|x|^{-bp})} 
\|\xi^2 u\|^2_{\mathcal{D}^{1,2}_{a}(\Omega)}}
{\|\nabla u \xi\|^2_{L^{2}(B_r(x_0),\di\mu_a)}}
+ 4r^{-2}\|u\|^2_{L^{2}(B_r(x_0),\di\mu_a)}\\
&\le 
2{\mathcal C}^2_{a,b,N}\|f\|^2_{L^{\frac{p}{p-1}}(\Omega,|x|^{-bp})} 
\left( \frac{8r^{-2}\|u\|^2_{L^{2}(B_r(x_0),\di\mu_a)}}
{\|\nabla u \xi\|^2_{L^{2}(B_r(x_0),\di\mu_a)}}+2\right)\\
&\quad+ 4r^{-2}\|u\|^2_{L^{2}(B_r(x_0),\di\mu_a)}.
\end{align*}
Thus
\begin{align*}
\int_{B_r(x_0)} |\nabla u|^2 \xi^2 \di \mu_a 
&\le 8{\mathcal C}^2_{a,b,N}\|f\|^2_{L^{\frac{p}{p-1}}(\Omega,|x|^{-bp})}
+ 4r^{-2}\|u\|^2_{L^{2}(B_r(x_0),\di\mu_a)}.
\end{align*}
Taking $r=r_0$ we have for $0<\rho\le r_0/2$
\begin{align*}
\media_{B_{\rho}(x_0)} |\n u|^2\di\mu_a\leq
C(N,a,\Omega,r_0) \bigg(\int_{\Omega}|u|^2\di\mu_a + 
\|f\|^2_{L^{\frac{2p}{p-2}}(\Omega,|x|^{-bp}\di x)}\bigg) \rho^{-2+2\a}\\
\end{align*}
From the above estimate, Corollary \ref{c:avera} and, the fact that $(N-bp)2/p=N-2-2a$ 
we derive the desired conclusion.
\end{proof}

\section{A Brezis-Kato type Lemma}
As in \cite{BrezisKato} we prove the following lemma to start an iteration procedure. 
\begin{lemma}\label{t:breziskato}
Let $\Omega\subset \rz^N$ be open, $a$, $b$ and $p$ satisfy (\ref{eq:24}), and $q>2$. 
Suppose $u\in D^{1,2}_a(\rz^N)\cap L^q(\Omega, |x|^{-bp})$ is a weak solution of
\begin{equation}
\label{eq:2bk}
-\div\big(|x|^{-2a}\n u\big)-
 \frac{V(x)}{|x|^{bp}}\,u=\frac{f(x)}{|x|^{bp}}\quad\text{in}\ \Omega,
\end{equation}
where $f\in L^q(\Omega, |x|^{-bp})$ and 
$V$ satisfies for some $\ell>0$  
\begin{align}\label{eq:20}
\int_{|V(x)|\geq \ell}|x|^{-bp}|V|^{\frac p{p-2}}+\int_{\Omega\setminus
  B_{\ell}(0)}|x|^{-bp}|V|^{\frac p{p-2}} \leq \min\bigg\{\frac18 {\mathcal
  C}_{a,b}^{-1},\frac{2}{q+4}{\mathcal C}_{a,b,N}^{-1} \bigg\}^{\frac p{p-2}}.
\end{align}
Then for any $\Omega' \Subset \Omega$
\begin{equation}\label{eq:breziskato}
\|u\|_{L^{\frac{pq}2}(\Omega', |x|^{-bp})}\leq C(\ell,q,\Omega')\|u\|_{L^{q}(\Omega,
  |x|^{-bp})}+\|f\|_{L^{q}(\Omega, |x|^{-bp})}.
\end{equation}
If, moreover, $u \in \Di$ then (\ref{eq:breziskato}) remains true for $\Omega'=\Omega$.
\end{lemma}
\begin{proof}
H\"older's inequality, (\ref{eq:25}) and (\ref{eq:20}) give for any $w\in\Di$  
\begin{align}
\notag
\int_{\Omega}|x|^{-bp}|V(x)|w^2 
&\le\ell\int\limits_{\substack{|V(x)|\leq\ell\ \text{and}\\ x\in\Omega\cap B_{\ell}(0)}} |x|^{-bp}w^2
+\int\limits_{\substack{|V(x)|\geq \ell\ \text{or}\\ x\in\Omega\setminus B_{\ell}(0)}} |x|^{-b(p-2)}|V||x|^{-2b}w^2\\
\notag
&\leq
\ell\int_{\Omega\cap B_{\ell}(0)}|x|^{-bp}w^2
+\bigg(\int_{\Omega}\frac{w^p}{|x|^{bp}}\bigg)^{\frac2p} \bigg(\int\limits_{\substack{|V(x)|\geq \ell\ \text{or}\\ x\in\Omega\setminus B_{\ell}(0)}} |x|^{-bp}|V|^{\frac{p}{p-2}}\bigg)^{\frac{p-2}p}\\
\label{eq:4bk}
&\leq \ell\int_{\Omega\cap B_{\ell}(0)}|x|^{-bp}w^2
+\min\Big(\frac18,\frac{2}{q+4}\Big)
\int_{\Omega}|x|^{-2a}|\n w|^2.
\end{align}
Suppose now that $u\in L^q(\Omega, |x|^{-bp})$. Fix $\Omega'\Subset \Omega$ and a nonnegative
cut-off function $\eta$, such that $\supp(\eta) \Subset \Omega$ and $\eta \equiv 1$ on $\Omega'$. 
Set $u^n := \min(n,|u|) \in D^{1,2}_a(\rz^N)$ and test (\ref{eq:2bk}) with $u(u^n)^{q-2}\eta^2
\in\Di$. This leads to 
\begin{align*}
(q&-2)\int_{\Omega}|x|^{-2a}\eta^2|\n u^n|^2(u^n)^{q-2}+ 
\int_{\Omega}|x|^{-2a}\eta^2 (u^n)^{q-2}|\n u|^2\\   
&=\int_{\Omega}|x|^{-bp}V(x)\eta^2 u^2(u^n)^{q-2}
+\int_{\Omega}|x|^{-bp}f\eta^2 (u^n)^{q-2}u
-2\int_{\Omega} |x|^{-2a} \n u \eta  (u_n)^{q-2} \n \eta u.
\end{align*}
We use the elementary inequality $2ab \le 1/2 a^2 + 4 b^2$ and obtain
\begin{align}
\label{eq:19}
(q&-2)\int_{\Omega}|x|^{-2a}\eta^2|\n u^n|^2(u^n)^{q-2}+\frac12 
\int_{\Omega}|x|^{-2a}\eta^2 (u^n)^{q-2}|\n u|^2\\
\notag   
&\le \int_{\Omega}|x|^{-bp}V(x)\eta^2 u^2(u^n)^{q-2}
+\int_{\Omega}|x|^{-bp}f\eta^2 (u^n)^{q-2}u
+ 4 \int_{\Omega} |x|^{-2a}|\n \eta|^2 u^2(u_n)^{q-2}.
\end{align}
Furthermore, an explicit calculation gives
\begin{align}\label{eq:21}
\begin{split}
\big|\n \big((u^n)^{\frac q2-1}u\eta)\big|^2
&\le \frac{(q+4)(q-2)}4(u^n)^{q-2}\eta^2|\n u^n|^2+2(u^n)^{q-2}
|\n u|^2\eta^2\\
&\quad+2 (u^n)^{q-2} u^2|\n \eta|^2+\frac{q-2}2 (u^n)^{q}|\n \eta|^2.
\end{split}
\end{align}
Let $C(q):=\min\big\{\frac14, \frac4{q+4}\big\}$. From (\ref{eq:19}) and (\ref{eq:21}) we get
\begin{align}\label{eq:22}
C(q)\int_{\Omega}
\frac{\big|\n \big((u^n)^{\frac q2-1}u\eta\big)\big|^2}
{|x|^{2a}}
& \leq 2(2+C(q))
\int_{\Omega}\frac{(u^n)^{q-2}u^2|\n \eta|^2}{|x|^{2a}} +C(q)\frac{q-2}2\int_{\Omega}\frac{(u^n)^q|\n \eta|^2}{|x|^{2a}}\notag\\
&\quad +\int_{\Omega} \frac {f(x)}{|x|^{bp}}\eta^2 (u^n)^{q-2}u+\int_{\Omega} \frac
{V(x)}{|x|^{bp}}\eta^2 (u^n)^{q-2}u^2.
\end{align}
Estimate (\ref{eq:4bk}) applied to $\eta(u^n)^{\frac q2-1}u$ gives 
\begin{align}
\label{eq:5}
\int_{\Omega} \frac{V^+\big(\eta(u^n)^{\frac q2-1}u\big)^2}{|x|^{bp}}
&\leq
\frac{C(q)}{2}\int_{\Omega}\frac{\big|\n(\eta(u^n)^{\frac q2-1}u)\big|^2}{|x|^{2a}}
+\ell\int_{\Omega\cap B_{\ell}(0)}\frac{(u^n)^{q-2}u^2\eta^2}{|x|^{bp}}. 
\end{align}
By H\"older's inequality and convexity we arrive at
\begin{align}
\label{eq:10}
\int_{\Omega}\frac{|f|\eta}{|x|^{\frac{bp}q}}\frac{(u^n)^{q-2}u\eta}{|x|^{\frac{bp(q-1)}q}}
\leq 
\frac{q-1}q\int_{\Omega}|x|^{-bp}\eta^{\frac q{q-1}}|u^n|^{q\frac{q-2}{q-1}}|u|^{\frac q{q-1}} 
+\frac1q\int_{\Omega}|x|^{-bp}|f|^q\eta^q.
\end{align}
We use (\ref{eq:5}) and (\ref{eq:10}) to estimate the terms with $f$ and $V$ in (\ref{eq:22}), 
then (\ref{eq:25}) yields
\begin{align*}
\bigg(\int_{\Omega}&|x|^{-bp}|u^n|^{(\frac q2-1)p}|u|^p\eta^p\bigg)^{\frac 2p}\\
&\le \frac{2{\mathcal C}_{a,b,N}(q-1)}{C(q)q}\int_{\Omega}|x|^{-bp}\eta^{\frac
  q{q-1}}|u^n|^{q\frac{q-2}{q-1}}|u|^{\frac q{q-1}}+\frac{2{\mathcal
    C}_{a,b}}{C(q)q}\int_{\Omega}|x|^{-bp}|f(x)|^q\eta^q\\
&\quad +
\frac{2\ell{\mathcal C}_{a,b,N}}{C(q)}\int_{\Omega\cap
  B_{\ell}(0)}|x|^{-bp}\eta^2|u^n|^{q-2}u^2+\frac{4{\mathcal
    C}_{a,b}(2+C(q))}{C(q)}\int_{\Omega}|x|^{-2a}|u^n|^{q-2}u^2|\n \eta|^2\\
&\quad +{\mathcal C}_{a,b,N}(q-2)\int_{\Omega}|x|^{-2a} |u^n|^q |\n \eta|^2.
\end{align*}
Letting $n\to\infty$ in the above inequality (\ref{eq:breziskato}) follows. 
Observe that if $u \in \Di$ then we need not to use the cut-off
$\eta$ and the same analysis as above gives the estimate (\ref{eq:breziskato}) for $\Omega'=\Omega$.
The lemma is thereby proved. 
\end{proof}
\begin{remark}
By Vitali's theorem $V$ belongs to $L^{p/(p-2)}(\Omega,|x|^{-bp})$ 
if and only if there exists $\ell$ such that (\ref{eq:20}) is satisfied. But the constant in (\ref{eq:breziskato}) 
depends uniformly on $\ell$ and not on the norm of $V$ in $L^{p/(p-2)}(\Omega,|x|^{-bp})$.  
\end{remark}

\begin{proof}[{\bf Proof of Theorem \ref{t:reg}}]
We apply Lemma \ref{t:breziskato} with $f=0$ and
$V(x)=K(x)|u|^{p-2}$. Starting with $q=p$, the lemma gives $u\in
L_{\rm loc}^{\frac{p^2}{2}}(\Omega,|x|^{-bp})$. Taking $q=\frac{p^2}2$, we find $u\in
L_{\rm loc}^{\frac{p^3}{4}}(\Omega,|x|^{-bp})$. Iterating the process, we obtain
that $u \in
L_{\rm loc}^{p^{k+1}/2^k}(\Omega,|x|^{-bp})$ for any~$k$. Let $k_0\in\N$ be such that
$(p/2)^{k_0}\geq 2(p-1)/(p-2)$, then after $k_0$ steps we find that $u\in L_{\rm loc}
  ^{\frac{2p(p-1)}{p-2}}(\Omega)$. Having this high integrability we may use
Theorem \ref{t:holder_estimate} with $f(x)=K(x) |u|^{p-2}u$ to get the desired regularity of $u$.  
\end{proof}

\appendix

\section{}

\begin{lemma}\label{l:veronica}
Let $a$, $b$ and $p$ satisfy (\ref{eq:24}) and $\eps>0$. 
Then we have
\begin{align}
\label{eq:9}
\left(\int_{B_{\rho}(x_0)}|x|^{-bp}\right)^{2/p+\eps}\leq C_{(\ref{eq:9})}(N)\,
\rho^{-2+\eps N} (\max(\rho,|x_0|)^{-\eps bp} \int_{B_{\rho}(x_0)}|x|^{-2a}.    
\end{align}
\end{lemma}
\begin{proof}
Let us distinguish two cases.\\
{\bf Case 1:} $\rho\geq |x_0|/2$. Since 
$(N-bp)({2}/{p}+\eps)= N-2-2a+\eps(N-bp)$
we obtain
\begin{align*}
\left(\int_{B_{\rho}(x_0)}|x|^{-bp}\right)^{2/p+\eps}
&\leq\left(\int_{B_{3\rho}(0)}|x|^{-bp}\right)^{2/p+\eps}=C_1(N)\rho^{N-2-2a+\eps(N-bp)}.
\end{align*}
From the doubling property (\ref{eq:3}) and the
fact that $B_{\rho}(0)\subset B_{4\rho}(x_0)$ we infer,
\begin{align*}
\rho^{\eps(N-bp)-2}\int_{B_{\rho}(x_0)}|x|^{-2a}&\geq
c\rho^{\eps(N-bp)-2}\int_{B_{4\rho}(x_0)}|x|^{-2a}\\
&\geq c\rho^{\eps(N-bp)-2}\int_{B_{\rho}(0)}|x|^{-2a}=C_2(N) \rho^{N-2-2a+\eps(N-bp)}
\end{align*}
and the claim follows in Case 1.\\
{\bf Case 2:} $\rho< |x_0|/2$. We have for all $x \in B_r(x_0)$ that
${1}/{2} \le {|x|}/{|x_0|} \le 3$.
Consequently,
\begin{align*}
\left(\int_{B_{\rho}(x_0)}|x|^{-bp}\right)^{2/p+\eps}
&\leq C_3(N)\rho^{N(2/p+\eps)} |x_0|^{-2b-\eps bp}\\
&\le C_3(N)\rho^{N-2} |x_0|^{-2a} \rho^{2N/p-N+2} |x_0|^{-2(b-a)} \rho^{N\eps} 
|x_0|^{-\eps bp} 
\end{align*}
From $r<|x_0|/2$ we get
\begin{align*}
\left(\int_{B_{\rho}(x_0)}|x|^{-bp}\right)^{2/p+\eps}
&\le C_3(N)\rho^{N-2} |x_0|^{-2a} \rho^{N\eps} |x_0|^{-\eps bp}\\
&\le C_4(N)  \left(\int_{B_{\rho}(x_0)}|x|^{-2a}\right)\rho^{N\eps} |x_0|^{-\eps bp},
\end{align*}
which ends the proof.
\end{proof}

\begin{lemma}
\label{reg:lem:1}
Suppose $\Phi$ be a nonnegative and nondecreasing functions on $[0,R]$ such that
\begin{align}
\label{eq:1}
\Phi(\rho) &\le A_1~ \mu_a\big(B_\rho(x)\big) \mu_a\big(B_r(x)\big)^{-1}
\big(\frac{\rho}{r}\big)^{-\alpha} \Phi(r) + A_2~ \mu_a\big(B_r(x)\big) r^{-\beta},    
\end{align}
for any $0<\rho\le r\le R$, where $A_1$, $A_2$, $\alpha$ and $\beta$ are positive constants
satisfying $\alpha<\beta$. Then for any $\gamma \in (\alpha,\beta)$ there exists a constant
$C_{(\ref{eq:2})}=C_{(\ref{eq:2})}(A_1,\alpha,\beta,\gamma)$ independent of $x$ and $r$ such that for $0<\rho\le
r \le R$
\begin{align}
\label{eq:2}
\Phi(\rho) &\le C_{(\ref{eq:2})} \left(\mu_a\big(B_\rho(x)\big) \mu_a\big(B_r(x)\big)^{-1}
\big(\frac{\rho}{r}\big)^{-\gamma} \Phi(r) + A_2~ \mu_a\big(B_\rho(x)\big) \rho^{-\beta}\right).    
\end{align}
\end{lemma}
\begin{proof}
Fix $\gamma \in (\alpha,\beta)$ and set $\tau:= \min(A_1^{-1/(\gamma-\alpha)},1/2)$.
Then we have for $0<r\le R$
\begin{align*}
\Phi(\tau r) \le \mu_a(B_{\tau r}(x)) \mu_a(B_r(x))^{-1} \tau^{-\gamma}
\Phi(r) + A_2 r^{-\beta} \mu_a(B_r(x)).      
\end{align*}
Hence we may estimate for $k\in \nz$
\begin{align*}
\Phi(\tau^{k+1} r) &\le  
\mu_a(B_{\tau^{k-1} r}(x)) \mu_a(B_{\tau^{k} r}(x))^{-1} \tau^{\beta-\gamma}\big)\\
&\le  \mu_a(B_{\tau^{k+1} r}(x)) \mu_a(B_r(x))^{-1} \tau^{-(k+1)\gamma} \Phi(r) 
+ A_2(\tau^{k} r)^{-\beta} \mu_a(B_{\tau^{k} r}(x))\\
&\quad \cdot \sum_{j=0}^k  \underbrace{\mu_a(B_{\tau^{k+1} r}(x)) \mu_a(B_{\tau^{k} r}(x))^{-1}}_{\le 1}\, 
\underbrace{\mu_a(B_{\tau^{k-j} r}(x)) \mu_a(B_{\tau^{k-j+1} r}(x))^{-1}}_{\le C_{(\ref{eq:3})}(\tau) \text{ by
    (\ref{eq:3})}} \tau^{(\beta-\gamma)j} \\
&\le C_{(\ref{eq:3})}(\tau)\mu_a(B_{\tau^{k+2} r}(x)) \mu_a(B_r(x))^{-1} \tau^{-(k+1)\gamma} \Phi(r) 
+  \frac{A_2 C_{(\ref{eq:3})}(\tau)}{1-\tau^{\beta-\gamma}} (\tau^{k} r)^{-\beta} \mu_a(B_{\tau^{k} r}(x))   
\end{align*}
For $0<\rho\le r$ we may choose $k \in \nz$ such that $\tau^{k+2}r<\rho<\tau^{k+1}r$ and obtain
\begin{align*}
\Phi(\rho) &\le  \Phi(\tau^{k+1} r)\\
&\le C_{(\ref{eq:3})}(\tau)\mu_a(B_\rho(x)) \mu_a(B_r(x))^{-1} \big(\frac{\rho}{r}\big)^{-\gamma} \Phi(r) 
+ \frac{A_2 C_{(\ref{eq:3})}^3(\tau)}{\tau (1-\tau^{\beta-\gamma})} \mu_a(B_\rho(x)) \rho^{-\beta}.
\end{align*}
\end{proof}

\bibliographystyle{adinat}
\bibliography{regularity}

\def\cprime{$'$}
\begin{thebibliography}{8}
\expandafter\ifx\csname natexlab\endcsname\relax\def\natexlab#1{#1}\fi
\expandafter\ifx\csname url\endcsname\relax
  \def\url#1{{\tt #1}}\fi
\expandafter\ifx\csname urlprefix\endcsname\relax\def\urlprefix{URL }\fi

\bibitem[{Br{\'e}zis and Kato(1979)}]{BrezisKato}
H.~Br{\'e}zis and T.~Kato.
\newblock {\em Remarks on the {S}chr\"odinger operator with singular complex
  potentials\/}.
\newblock J. Math. Pures Appl. (9), {\bf 58} (1979), no.~2, 137--151.

\bibitem[{Caffarelli et~al.(1984)Caffarelli, Kohn and Nirenberg}]{CKN}
L.~Caffarelli, R.~Kohn and L.~Nirenberg.
\newblock {\em First order interpolation inequalities with weights\/}.
\newblock Compositio Math., {\bf 53} (1984), no.~3, 259--275.

\bibitem[{Catrina and Wang(2001)}]{CatrinaWang}
F.~Catrina and Z.-Q. Wang.
\newblock {\em On the {C}affarelli-{K}ohn-{N}irenberg inequalities: sharp
  constants, existence (and nonexistence), and symmetry of extremal
  functions\/}.
\newblock Comm. Pure Appl. Math., {\bf 54} (2001), no.~2, 229--258.

\bibitem[{Colorado and Peral(2003)}]{ColoPeral}
E.~Colorado and I.~Peral.
\newblock {\em Eigenvalues and bifurcation for elliptic equations with mixed
  {D}irichlet-{N}eumann boundary conditions related to
  {C}affarelli-{K}ohn-{N}irenberg inequalities\/}.
\newblock  (2003).
\newblock Preprint.

\bibitem[{De~Cicco and Vivaldi(1996)}]{CiccoVivaldi96}
V.~De~Cicco and M.~A. Vivaldi.
\newblock {\em Harnack inequalities for {F}uchsian type weighted elliptic
  equations\/}.
\newblock Comm. Partial Differential Equations, {\bf 21} (1996), no. 9-10,
  1321--1347.

\bibitem[{Fabes et~al.(1982)Fabes, Kenig and Serapioni}]{FabesKenigSera82}
E.~B. Fabes, C.~E. Kenig and R.~P. Serapioni.
\newblock {\em The local regularity of solutions of degenerate elliptic
  equations\/}.
\newblock Comm. Partial Differential Equations, {\bf 7} (1982), no.~1, 77--116.

\bibitem[{Felli and Schneider(2003)}]{SF_bup}
V.~Felli and M.~Schneider.
\newblock {\em Compactness and existence results for degenerate critical
  elliptic equations\/} (2003).
\newblock Preprint.

\bibitem[{Heinonen et~al.(1993)Heinonen, Kilpel{\"a}inen and
  Martio}]{potential}
J.~Heinonen, T.~Kilpel{\"a}inen and O.~Martio.
\newblock {\em Nonlinear potential theory of degenerate elliptic equations\/}.
\newblock The Clarendon Press Oxford University Press, New York (1993).

\end{thebibliography}

\end{document}